\documentclass[12pt,reqno]{amsart}
\usepackage{amscd,amsmath,amsthm,amssymb}
\usepackage{color}
\usepackage{pstricks}
\usepackage{stmaryrd}
\usepackage{tikz}
\usepackage{url}
\usepackage{bm}

\usepackage{latexsym}
\usepackage{amsfonts,amsmath,mathtools}
\usepackage{graphics}
\usepackage{float}
\usepackage{enumitem}

\usepackage{booktabs} 
\usepackage{colortbl}
\usepackage{xcolor}

\newpsstyle{fatline}{linewidth=1.5pt}
\newpsstyle{fyp}{fillstyle=solid,fillcolor=verylight}
\definecolor{verylight}{gray}{0.97}
\definecolor{light}{gray}{0.9}
\definecolor{medium}{gray}{0.85}
\definecolor{dark}{gray}{0.6}

\def\NZQ{\mathbb}               

\def\FF{{\NZQ F}}

\def\G{{\mathcal G}}

\def\HS{\textup{HS}}
\def\pd{\textup{proj}\phantom{.}\!\textup{dim}}

\def\opn#1#2{\def#1{\operatorname{#2}}} 
\opn\chara{char} \opn\length{\ell} \opn\pd{pd} \opn\rk{rk}
\opn\projdim{proj\,dim} \opn\injdim{inj\,dim} \opn\rank{rank}
\opn\depth{depth} \opn\grade{grade} \opn\height{height}
\opn\embdim{emb\,dim} \opn\codim{codim}

\opn\Tr{Tr} \opn\bigrank{big\,rank}
\opn\superheight{superheight}\opn\lcm{lcm}
\opn\trdeg{tr\,deg}
\opn\reg{reg} \opn\lreg{lreg} \opn\ini{in} \opn\lpd{lpd}
\opn\size{size} \opn\sdepth{sdepth}
\opn\link{link}\opn\fdepth{fdepth}\opn\lex{lex}
\opn\tr{tr}
\opn\type{type}
\opn\gap{gap}
\opn\diam{diam}
\opn\Mod{Mod}
\opn\div{div} \opn\Div{Div} \opn\cl{cl} \opn\Cl{Cl}
\opn\Spec{Spec} \opn\Supp{Supp} \opn\supp{supp} \opn\Sing{Sing}
\opn\Ass{Ass} \opn\Min{Min}\opn\Mon{Mon}
\opn\Ann{Ann} \opn\Rad{Rad} \opn\Soc{Soc}
\opn\Im{Im} \opn\Ker{Ker} \opn\Coker{Coker} \opn\Am{Am}
\opn\Hom{Hom} \opn\Tor{Tor} \opn\Ext{Ext} \opn\End{End}
\opn\Aut{Aut} \opn\id{id}

\opn\nat{nat}
\opn\pff{pf}
\opn\Pf{Pf} \opn\GL{GL} \opn\SL{SL} \opn\mod{mod} \opn\ord{ord}
\opn\Gin{Gin} \opn\Hilb{Hilb}\opn\sort{sort}
\opn\PF{PF}\opn\Ap{Ap}
\opn\dist{dist}
\opn\aff{aff}
\opn\relint{relint} \opn\st{st}
\opn\lk{lk} \opn\cn{cn} \opn\core{core} \opn\vol{vol}  \opn\inp{inp} \opn\nilpot{nilpot}
\opn\link{link} \opn\star{star}\opn\lex{lex}\opn\set{set}
\opn\width{wd}
\opn\Fr{F}
\opn\QF{QF}
\opn\G{G}
\opn\type{type}\opn\res{res}
\opn\conv{conv}
\opn\sr{sr}
\opn\gr{gr}

\def\pot#1#2{#1[\kern-0.28ex[#2]\kern-0.28ex]}

\opn\dirlim{\underrightarrow{\lim}}
\opn\inivlim{\underleftarrow{\lim}}
\def\Implies{\ifmmode\Longrightarrow \else
	\unskip${}\Longrightarrow{}$\ignorespaces\fi}
\def\implies{\ifmmode\Rightarrow \else
	\unskip${}\Rightarrow{}$\ignorespaces\fi}
\def\iff{\ifmmode\Longleftrightarrow \else
	\unskip${}\Longleftrightarrow{}$\ignorespaces\fi}

\let\:=\colon
\newtheorem{Theorem}{Theorem}[section]
\newtheorem{Lemma}[Theorem]{Lemma}
\newtheorem{Corollary}[Theorem]{Corollary}
\newtheorem{Proposition}[Theorem]{Proposition}

\newtheorem{Example}[Theorem]{Example}

\newtheorem{Definition}[Theorem]{Definition}

\newtheorem{Conjecture}[Theorem]{Conjecture}

\let\epsilon\varepsilon
\let\kappa=\varkappa
\def\qed{\ifhmode\textqed\fi
	\ifmmode\ifinner\hfill\quad\qedsymbol\else\dispqed\fi\fi}
\def\textqed{\unskip\nobreak\penalty50
	\hskip2em\hbox{}\nobreak\hfill\qedsymbol
	\parfillskip=0pt \finalhyphendemerits=0}
\def\dispqed{\rlap{\qquad\qedsymbol}}

\opn\dis{dis}
\def\pnt{{\raise0.5mm\hbox{\large\bf.}}}

\opn\Lex{Lex}
\opn\soc{soc}
\def\set{\textup{set}}

\usepackage{lipsum}

\begin{document}

	\title{Homological Shifts of Polymatroidal ideals}
	\author{Antonino Ficarra}
	
	\address{Antonino Ficarra, Department of mathematics and computer sciences, physics and earth sciences, University of Messina, Viale Ferdinando Stagno d'Alcontres 31, 98166 Messina, Italy}
	\email{antficarra@unime.it}
	
	\thanks{.
	}
	
	
	
	\maketitle
	
	\begin{abstract}
		We study the homological shifts of polymatroidal ideals. In our main theorem we prove that the first homological shift ideal of any polymatroidal ideal is again polymatroidal, supporting a conjecture of Bandari, Bayati and Herzog that predicts that all homological shift ideals of a polymatroidal ideal are polymatroidal. As a nice consequence, we recover a result of Bayati which proves this conjecture in the squarefree case.
	\end{abstract}
	
	\section{Introduction}
	The study of the minimal free resolutions of monomial ideals is a main topic in combinatorial commutative algebra. Let $K$ be a field, $I$ a monomial ideal of the standard graded polynomial ring $S=K[x_1,\dots,x_n]$ and $\FF$ be its minimal multigraded free $S$-resolution. Then, the free modules in the resolution $\FF$ are of the form $F_i=\bigoplus_{j=1}^{\beta_i(I)}S(-{\bf a}_{i,j})$, where ${\bf a}_{i,j}$ are integral vectors with non negative entries, called the \textit{multigraded shifts} of $I$. For an integral vector ${\bf a}=(a_1,a_2,\dots,a_n)$, we let ${\bf x^a}=x_1^{a_1}x_{2}^{a_2}\cdots x_{n}^{a_n}$. Using this notation, the \textit{$i$th homological shift ideal} of $I$ is defined as $\HS_i(I)=({\bf x}^{{\bf a}_{i,j}}:j=1,\dots,\beta_i(I))$, see \cite{HMRZ021a}. Note that $\HS_0(I)=I$ and $\HS_j(I)=(0)$ for $j>\pd(I)$.
	
	The basic goal of this theory is to understand how much of the homological and combinatorial properties of a monomial ideal $I$ are inherited by its homological shift ideals. We refer to any combinatorial or homological property enjoyed by $\HS_0(I)=I$ and by $\HS_j(I)$, for all $1\le j\le\pd(I)$, as an \textit{homological shift property} of $I$. This rather new trend of research, see \cite{Bay019,Bay2023,BJT019,CF,CF1,F2Pack,FH2023,HMRZ021a,HMRZ021b}, has its origins in a meeting between Somayeh Bandari, Shamila Bayati and J\"urgen Herzog that took place in Essen in 2012. Due to experimental evidence, the three authors conjectured that the property of being \textit{polymatroidal} is an homological shift property. For a monomial $u$, we define the \textit{$x_i$-degree} of $u$ as $\deg_{x_i}(u)=\max\{j:x_i^j\ \text{divides}\ u\}$. Recall that a \textit{polymatroidal} ideal $I$ is an equigenerated monomial ideal of $S$ whose minimal generating set $G(I)$ corresponds to the base of a \textit{discrete polymatroid}. The bases of a discrete polymatroidal can be characterized in term of the so--called \textit{exchange property}. Thus an ideal $I\subset S$ is polymatroidal if and only if the following exchange property holds: for all $u,v\in G(I)$ such that $\deg_{x_i}(u)>\deg_{x_i}(v)$, there exists $j$ such that $\deg_{x_j}(u)<\deg_{x_j}(v)$ and $x_j(u/x_i)\in G(I)$. The Bandari--Bayati--Herzog conjecture is widely open. Until now it is solved only in two cases: for squarefree polymatroidal, \textit{i.e.}, \textit{matroidal ideals} \cite{Bay019}, and for polymatroidal ideals that satisfy the \textit{strong exchange property} \cite{HMRZ021a}.
	
	Polymatroidal ideals are
	one of the most distinguished classes of monomial ideals. Indeed, the product of polymatroidal ideals is polymatroidal. Any polymatroidal ideal has linear quotients and thus a linear resolution. Thus, they have \textit{linear powers}, \cite[Corollary 12.6.4]{JT}, a rare property among monomial ideals. In this paper our aim is to study the homological shift ideals of polymatroidal ideals. See also \cite{BH2013,BanRam019,HH2003,HHV2005,HRV2013,Lu2014,Schweig2011}.
	
	The paper is organized as follows. In Section \ref{sec:1}, we collect some basic facts on homological shift ideals of equigenerated monomial ideals with linear quotients. Recall that a monomial ideal $I\subset S$ has \textit{linear quotients} if for some \textit{admissible order} $u_1>\dots >u_m$ of its minimal generating set $G(I)$, the colon ideals $(u_1,\dots,u_{k-1}):u_k$ are generated by a subset of the variables, for $k=2,\dots,m$. Quite generally, to determine the ideals $\HS_j(I)$ for a monomial ideal $I$ is difficult. Nonetheless, this is easy for ideals with linear quotients as shown in Proposition \ref{Prop:MultidegreesIlinQuotEquig}.
	
	In Section \ref{sec:2}, we prove our main theorem. We are able to show that $\HS_1(I)$ is polymatroidal if $I$ is polymatroidal (Theorem \ref{Thm:HS1(I)PolyMatr}). This is the first result valid for \textit{all} polymatroidal ideals, that supports Conjecture \ref{Conj:BBH}. Our proof is based on Proposition \ref{prop:HS1(I)linquot}, whose main advantage consists in the fact that $\HS_{1}(I)$ does not depend upon the admissible order of $I$. To study the higher homological shift ideals, firstly we note that for $j\ge1$, $\HS_{j+1}(I)\subseteq(\HS_1(\HS_j(I)))_{>j+1}$, (Corollary \ref{cor:HS1HSjLexEquig}), where $J_{>j+1}$ is the monomial ideal with $G(J_{>j+1})=\{u\in G(J):|\supp(u)|>j+1\}$ as a minimal generating set and $\supp(u)=\{i:x_i\ \textup{divides}\ u\}$. Unfortunately, equality in the above inclusion does not hold in general, (Example \ref{Ex:HS1(HSj)notTransversal}).  Nonetheless, it holds for matroidal ideals, (Proposition \ref{prop:HSj+1=HS1HSjmatroidal}). As a consequence Conjecture \ref{Conj:BBH} holds for all matroidal ideals, (Corollary \ref{cor:HSj(I)matroidal}). It would be of interest to classify all polymatroidal ideals satisfying the equation $\HS_{j+1}(I)=(\HS_1(\HS_j(I)))_{>j+1}$ for all $j<\pd(I)$.
	
	\section{Generalities on homological shift ideals}\label{sec:1}
	
	Let $S=K[x_1,\dots,x_n]$ be the standard graded polynomial ring over a field $K$. For a monomial $u=x_1^{a_1}\cdots x_{n}^{a_n}\in S$, the vector ${\bf a}=(a_1,\dots,a_n)$ is called the \textit{multidegree} of $u$. We put $u={\bf x^a}$. For ${\bf a}={\bf 0}=(0,0,\dots,0)$, ${\bf x^0}=1$. Whereas $\deg(u)=a_1+a_2+\dots+a_n$ is the \textit{degree} of $u$. Let $G(I)$ be the unique minimal monomial generating set of $I$, and $G(I)_d=\{u\in G(I):\deg(u)=d\}$.
	\begin{Definition}
		\rm Let $I\subset S$ be a monomial ideal, and let $(\mathbb{F},\partial)$ be the minimal multigraded free resolution of $I$. The \textit{$i$th homological shift ideal} of $I$ is defined as
		$$
		\HS_i(I)\ =\ ({\bf x^a}\ :\ \beta_{i,{\bf a}}(I)\ne0).
		$$
		Here $\beta_{i,{\bf a}}(I)$ is a multigraded Betti number.
	\end{Definition}
	
	Clearly $\HS_0(I)=I$ and $\HS_j(I)=(0)$ for $j>\pd(I)$. In general, we ask what properties of $I$ are inherited by its homological shift ideals.
	
	Let $I$ be a monomial ideal with minimal generating set $G(I)=\{u_1,\dots,u_m\}$. We say that $I$ has \textit{linear quotients} with respect to the order $u_1>\dots>u_m$ of its minimal generators, if for all $k=2,\dots,m$, the colon ideal $(u_1,\dots,u_{k-1}):u_k$ is generated by a subset of the variables, $x_1,\dots,x_n$. Any such ordering of $G(I)$ is called an \textit{admissible order} of $I$. Set
	$$
	\set(u_k)\ =\ \{i:x_i\in(u_1,\dots,u_{k-1}):u_k\},
	$$
	for $k=2,\dots,m$ and $\set(u_1)=\emptyset$. Furthermore, put ${\bf x}_F=\prod_{i\in F}x_i$ if $F$ is non empty, and ${\bf x}_\emptyset=1$. By \cite[Lemma 1.5]{ET}, we have the following result.
	\begin{Proposition}\label{Prop:MultidegreesIlinQuotEquig}
		Let $I\subset S$ be a monomial ideal with linear quotients with admissible order $u_1>\dots>u_m$ of $G(I)$. Then,
		\begin{equation}\label{eq:HSi(I)linquot}
			\HS_j(I)\ =\ ({\bf x}_Fu\ :\ u\in G(I),\ F\subseteq\set(u),\ |F|=j).
		\end{equation}
	\end{Proposition}
	
	Let $u={\bf x^a}\in S$, the \textit{$x_i$-degree} of $u$ is $\deg_{x_i}(u)=\max\{j:x_i^j\ \textup{divides}\ u\}=a_i$. Given monomials of the same degree $u,v\in S$, the \textit{distance} between $u$ and $v$ is
	$$
	d(u,v)\ =\ \frac{1}{2}\sum\limits_{i=1}^{n}\big|\deg_{x_i}(u)-\deg_{x_i}(v)\big|.
	$$
	\begin{Lemma}\label{Lem:d(u,v)=1}
		Let $u,v$ monomials of $S$ of the same degree. Then, $u=x_k(v/x_{\ell})$ for some $k\ne\ell$, if and only if $d(u,v)=1$.
	\end{Lemma}
	\begin{proof}
		Let $v=x_{1}^{b_1}\cdots x_n^{b_n}$, then $u=x_k(v/x_{\ell})=(\prod_{i\ne k,\ell}x_i^{b_i})x_k^{b_k+1}x_{\ell}^{b_\ell-1}$. Note that
		$$
		\big|\deg_{x_i}(u)-\deg_{x_i}(v)\big|\ =\ \begin{cases}
			|b_{\ell}-1-b_{\ell}|&\text{if}\ i=\ell,\\
			\hfil|b_{k}\!+1-b_{k}|&\text{if}\ i=k,\\
			\hfill 0&\text{otherwise}.
		\end{cases}
		$$
		Thus $d(u,v)=\frac{1}{2}\sum_{i}|\deg_{x_i}(u)-\deg_{x_i}(v)|=\frac{1}{2}\big(|b_k+1-b_k|+|b_{\ell}-1-b_{\ell}|\big)=\frac{1}{2}(1+1)=1$.\smallskip
		
		Conversely, assume $d(u,v)=1$. From the definition of $d(u,v)$ it is clear that either $u=x_k^2v$ or $u=x_k(v/x_\ell)$, for some $k\ne\ell$. But the first possibility does not occur, lest $\deg(u)>\deg(v)$. Therefore, the desired conclusion follows.
	\end{proof}
	
	Combining Lemma \ref{Lem:d(u,v)=1} with \cite[Proposition 1.3]{HMRZ021a} we have
	\begin{Proposition}\label{prop:HS1(I)linquot}
		Let $I\subset S$ be an equigenerated monomial ideal with linear quotients. Then $$\HS_1(I)\ =\ (\lcm(u,v)\ :\ u,v\in G(I),\ d(u,v)=1).$$
	\end{Proposition}
	
	Let $\succ$ be a \textit{monomial order} on $S$. Up to a relabeling on the variables, we may assume that $x_1\succ\dots\succ x_n$. We say that $\succ$ \textit{is induced by} $x_1>\dots>x_n$. We say that $I$ has linear quotients with respect to $\succ$ if $I$ has linear quotients with admissible order $u_1\succ u_2\succ\dots\succ u_m$ of $G(I)$. A particular monomial order is the \textit{lex order} $>_{\lex}$ induced by $x_1>\dots>x_n$. Let ${\bf x^a},{\bf x^b}$ be monomials of $S$. Then ${\bf x^a}>_{\lex}{\bf x^b}$ if $a_1=b_1$, $\dots$, $a_{s-1}=b_{s-1}$ and $a_s>b_s$, for some $1\le s\le n$.
	
	We set $[n]=\{1,2,\dots,n\}$. For a monomial $u\in S$, we define its \textit{support} as $\supp(u)=\{i\in[n]:\deg_{x_i}(u)>0\}$ and we set $\max(u)=\max\supp(u)$.
	\begin{Lemma}\label{lem:set(u)>lexLinQuot}
		Let $I\subset S$ be an equigenerated monomial ideal with linear quotients with respect to $\succ$ induced by $x_1>\dots>x_n$. Then, for all $u\in G(I)$,
		$$
		\set(u)\ \subseteq\ [\max(u)-1].
		$$
	\end{Lemma}
	\begin{proof}
		Indeed, let $G(I)$ ordered as $u_1\succ\cdots\succ u_m$ and let $j\in\{1,\dots,m\}$. If $i\in\set(u_j)$, then $x_iu_j\in(u_1,\dots,u_{j-1})$. Since $\deg(u_1)=\dots=\deg(u_{j-1})=\deg(u_j)$, there exists $s\in\supp(u_j)$, $s\ne i$ such that $x_i(u_j/x_s)=u_p$ for some $p\le j-1$. But $x_i(u_j/x_s)=u_p\succ u_j$. Since $\succ$ is a monomial order, this implies that $x_i u_j\succ x_su_j$, thus $x_i\succ x_s$. Hence $i<s$. But $s\le\max(u_j)$ and so $i<\max(u_j)$.
	\end{proof}
	
	In \cite{HMRZ021a}, the following general inclusion was shown.
	\begin{Proposition}\label{Prop:InclusionHSLinQuot}
		\textup{\cite[Proposition 1.4]{HMRZ021a}} Let $I\subset S$ be a monomial ideal with linear quotients. Then $\HS_{j+1}(I)\subseteq\HS_1(\HS_j(I))$, for all $j$.
	\end{Proposition}
	
	However, in general $\HS_{j+1}(I)\ne\HS_1(\HS_j(I))$. Let $I=(x_2x_4,x_1x_2,x_1x_3)$. Then $I$ has linear quotients with admissible order $x_2x_4>x_1x_2>x_1x_3$. We have $\HS_1(I)=(x_1x_2x_3,x_1x_2x_4)$ and $\HS_2(I)=(0)$, but $\HS_1(\HS_1(I))=(x_1x_2x_3x_4)\ne(0)=\HS_2(I)$.\medskip
	
	For a monomial ideal $J\subset S$, let $J_{>\ell}$ be the monomial ideal whose minimal generating set is $G(J_{>\ell})=\{u\in G(J):|\supp(u)|>\ell\}$.
	\begin{Corollary}\label{cor:HS1HSjLexEquig}
		Let $I\subset S$ be an equigenerated monomial ideal with linear quotients with respect to a monomial order $\succ$ \textup{(}e.g., $>_{\lex}$\textup{)} induced by $x_1>\dots>x_n$. Then,
		$$
		\HS_{j+1}(I)\ \subseteq\ \big(\HS_1(\HS_j(I))\big)_{>j+1}.
		$$
	\end{Corollary}
	\begin{proof}
		For $j=0$ the assertion is immediate. Let $j>0$. Firstly we show the inclusion $\HS_{j+1}(I)\subseteq\HS_1(\HS_j(I))$. By Proposition \ref{prop:HS1(I)linquot},
		$$
		\HS_1(\HS_j(I))=(\lcm(w_1,w_2)\ :\ w_1,w_2\in G(\HS_j(I)),\ d(w_1,w_2)=1).
		$$
		Take $w={\bf x}_Fu\in G(\HS_{j+1}(I))$ with $u\in G(I)$, $F\subseteq\set(u)$ and $|F|=j+1$. Since $j+1\ge2$ we can find $r,s\in F$, $r\ne s$. Then $w_1={\bf x}_{F\setminus\{r\}}u,w_2={\bf x}_{F\setminus\{s\}}u\in G(\HS_j(I))$ and $d(w_1,w_2)=1$. Thus $\lcm(w_1,w_2)=w\in G(\HS_1(\HS_j(I)))$, as desired. It remains to prove that any $w={\bf x}_Fu\in G(\HS_{j+1}(I))$ has $\supp(w)>j+1$. By Lemma \ref{lem:set(u)>lexLinQuot} we have $F\subseteq\set(u)\subseteq[\max(u)-1]$. Hence $\max(u)\notin F$ and $F\cup\{\max(u)\}\subseteq\supp(w)$, and since $|F|=j+1$ we have that $\supp(w)\ge|F|+1=j+2$, as desired.
	\end{proof}
	
	\section{The first homological shift ideal of polymatroidal ideals}\label{sec:2}
	
	Recall that an equigenerated monomial ideal $I\subset S$ is a \textit{polymatroidal} ideal if it satisfies the following \textit{exchange property},
	\begin{enumerate}
		\item[$(*)$] for all $u,v\in G(I)$ and all $i$ such that $\deg_{x_i}(u)>\deg_{x_i}(v)$, there exists $j$ with $\deg_{x_j}(u)<\deg_{x_j}(v)$ and such that $x_j(u/x_i)\in G(I)$.
	\end{enumerate}
	Such ideals are called \textit{polymatroidal} because their minimal generating set $G(I)$ corresponds to the basis of a \textit{discrete polymatroid}, see \cite[Chapter 12]{JT} for more information on this subject. For later use, we recall that a polymatroidal ideal satisfies also the following \textit{dual exchange property}, \cite[Lemma 2.1]{HH2003},\smallskip
	\begin{enumerate}
		\item[$(**)$] for all $u,v\in G(I)$ and all $j$ such that $\deg_{x_j}(u)<\deg_{x_j}(v)$, there exists $i$ with $\deg_{x_i}(u)>\deg_{x_i}(v)$ and such that $x_j(u/x_i)\in G(I)$.
	\end{enumerate}
	
	It is expected that the following is true.
	
	\begin{Conjecture}\label{Conj:BBH}
		\textup{(Bandari--Bayati--Herzog), \cite{Bay019,HMRZ021a}.} Let $I\subset S$ be a polymatroidal ideal. Then all homological shift ideals $\HS_{j}(I)$ are again polymatroidal, for all $j\ge0$.
	\end{Conjecture}
	
	It is known that a polymatroidal ideal $I\subset S$ has linear quotients with respect to the lex order induced by $x_1>\dots>x_n$, see \cite[Theorem 2.4]{BanRam019}. Using the description of $\HS_1(I)$ given in Proposition \ref{prop:HS1(I)linquot} we can prove
	
	\begin{Theorem}\label{Thm:HS1(I)PolyMatr}
		Let $I\subset S=K[x_1,\dots,x_n]$ be a polymatroidal ideal. Then $\HS_1(I)$ is a polymatroidal ideal.
	\end{Theorem}
	\begin{proof}
		By Proposition \ref{prop:HS1(I)linquot},
		\begin{align*}
			\HS_1(I)\ &=\ (x_iu\ :\ u\in G(I),\ i\in\set(u))\\&=\ (\lcm(u,v)\ :\ u,v\in G(I),\ d(u,v)=1).
		\end{align*}
		We must prove the following exchange property,\smallskip
		
		\begin{enumerate}
			\item[$(*)$] for all monomials $u,v\in G(I)$, all integers $k\in\set(u)$, $\ell\in\set(v)$ such that $u_1=x_{k}u\ne x_{\ell}v=v_1$ and all $i$ with $\deg_{x_i}(u_1)>\deg_{x_i}(v_1)$, there exists $j$ such that $\deg_{x_j}(u_1)<\deg_{x_j}(v_1)$ and $x_{j}(u_1/x_{i})=x_{j}(x_{k}u)/x_i\in G(\HS_1(I))$.
		\end{enumerate}
		
		We may assume that $i$ is different both from $k$ and $\ell$. Indeed, if $k=i$, then as $k\in\set(u)$ we have $u_1=x_{p}z$ for some $z\in G(I)\setminus\{u\}$, $p\ne k$, and we may use the element $x_{p}z$ with $p\ne k=i$. The same reasoning applies for $\ell$. In particular, since $i\ne k$, $i\ne\ell$ and by hypothesis $\deg_{x_i}(u_1)>\deg_{x_i}(v_1)$, we have $\deg_{x_i}(u)>\deg_{x_i}(v)$ as well. Since $I$ is polymatroidal, the set
		$$
		\Omega\ =\ \{h\in[n]\setminus\{i\}\ :\ \deg_{x_h}(u)<\deg_{x_h}(v)\ \text{and}\ x_{h}(u/x_{i})\in G(I)\}.
		$$
		is non empty. Let $h\in\Omega$ and set $w=x_{h}(u/x_{i})$. We distinguish two cases.
		
		\medskip\noindent
		\textsc{Case 1.} Suppose that $k\in\Omega$. For $h=k\in\Omega$, we have $\deg_{x_k}(u)<\deg_{x_k}(v)$ and $w=x_{k}(u/x_{i})\in G(I)$. We distinguish two more cases.
		
		\medskip\noindent
		\textsc{Subcase 1.1.} Assume $w=v$. By hypothesis $v_1=x_{\ell}v\in G(\HS_1(I))$. We show that the property $(*)$ is verified for the integer $j=\ell$. Indeed, as $i\ne\ell$ and $w=v$,
		\begin{align*}
			\deg_{x_\ell}(u_1)&=\deg_{x_\ell}(x_{k}u)=\deg_{x_\ell}(x_{k}u/x_{i})=\deg_{x_\ell}(w)\\&=\deg_{x_\ell}(v)<\deg_{x_\ell}(x_{\ell}v),
		\end{align*}
		and $v_1=x_{\ell}v=x_{\ell}w=x_{\ell}(x_{k}u)/x_{i}\in\ G(\HS_1(I))$, as desired.
		
		\medskip\noindent
		\textsc{Subcase 1.2.} Assume $w\ne v$. Thus, for some $r$, $\deg_{x_r}(w)>\deg_{x_r}(v)$. Since $I$ is polymatroidal, there exists an integer $m$ with $\deg_{x_m}(w)<\deg_{x_m}(v)$ and such that $w_1=x_{m}(w/x_{r})\in G(I)$. Clearly $m\ne r$. Hence, Lemma \ref{Lem:d(u,v)=1} implies that $d(w,w_1)=1$ and Proposition \ref{prop:HS1(I)linquot} implies that
		\begin{equation}\label{eq:subcase1.2HS1(Poly)}
			\lcm(w,w_1)=x_{m}w=x_{m}(x_{\ell}u)/x_{i}\in G(\HS_1(I)).
		\end{equation}
		It remains to prove that the integer $m$ satisfies the first condition of property $(*)$, namely $\deg_{x_m}(u_1)<\deg_{x_m}(v_1)$. First note that $m\ne i$. Lest, if $i=m$, by hypothesis $\deg_{x_i}(w)<\deg_{x_i}(v)$ and then $\deg_{x_i}(u)=\deg_{x_i}(w)+1\le\deg_{x_i}(v)$, against the fact that $\deg_{x_i}(u)>\deg_{x_i}(v)$. Thus, since $m\ne i$, we have
		$
		\deg_{x_m}(u)\le\deg_{x_m}(w)<\deg_{x_m}(v).
		$
		This inequality together with equation (\ref{eq:subcase1.2HS1(Poly)}) show that the integer $j=m$ satisfies the property $(*)$ in such a case.
		
		\medskip\noindent
		\textsc{Case 2.} Suppose that $k\notin\Omega$. Nonetheless, for some $h\in\Omega$, with $h\ne k$, we have $\deg_{x_h}(u)<\deg_{x_h}(v)$ and $w=x_{h}(u/x_{i})\in G(I)$. Since $k\in\set(u)$, there exist $z\in G(I)$, $z\ne u$ and $x_{p}\ne x_{k}$ such that $u_1=x_{k}u=x_{p}z$.
		
		\medskip\noindent
		\textsc{Subcase 2.1.} Suppose $d(w,z)=1$. As $h\in\Omega$ but $k\notin\Omega$ we have $h\ne k$. Thus, as $w=x_{h}(u/x_{i})$ and $z=x_{k}(u/x_{p})$ it follows that $p=i$, lest $d(w,z)>1$. Hence $p=i$ and Proposition \ref{prop:HS1(I)linquot} implies
		\begin{align*}
			\lcm(w,z)&=\lcm(x_{h}(u/x_{i}),x_{k}(u/x_{p}))=x_{h}(x_{k}u)/x_{i}\in G(\HS_1(I)).
		\end{align*}
		Finally, we just need to check that $\deg_{x_h}(u_1)<\deg_{x_h}(v_1)$. Indeed, as $h\ne k$,
		\begin{align*}
			\deg_{x_h}(x_{k}u)=\deg_{x_h}(u)<\deg_{x_h}(v)\le\deg_{x_h}(x_{\ell}v).
		\end{align*}
		
		\medskip\noindent
		\textsc{Subcase 2.2.} Suppose $d(w,z)>1$. Then $p\ne i$, lest $d(w,z)=1$ by \textsc{Subcase 2.1}. Thus $d(w,z)=d(x_h(u/x_i),x_k(u/x_p))=2$, $i\ne h$, $h\ne k$, $k\ne p$, $p\ne i$ and
		\begin{align*}
			\deg_{x_i}(w)&<\deg_{x_i}(z),& \deg_{x_h}(w)&>\deg_{x_h}(z),\\
			\deg_{x_k}(w)&<\deg_{x_k}(z),& \deg_{x_p}(w)&>\deg_{x_p}(z).
		\end{align*}
		Moreover, for all $q\ne i,h,k,p$ we have $\deg_{x_q}(w)=\deg_{x_q}(z)$. Since $w,z\in G(I)$ and $\deg_{x_i}(z)>\deg_{x_i}(w)$ we have $z_1=x_{h}(z/x_{i})\in G(I)$ or $z_2=x_{p}(z/x_{i})\in G(I)$. We distinguish two more cases.
		
		\medskip\noindent
		\textsc{Subcase 2.2.1.} Suppose $z_1=x_{h}(z/x_{i})\in G(I)$. Note that
		$$
		x_{p}(z_1/x_{k})=x_{p}(x_{h}(z/x_{i}))/x_{k}=x_{p}x_{h}x_{k}((u/x_{p})/x_{i})/x_{k}=x_{h}(u/x_i)=w.
		$$
		Since $k\ne p$, Lemma \ref{Lem:d(u,v)=1} implies that $d(z_1,w)=1$. Thus, by Proposition \ref{prop:HS1(I)linquot}
		\begin{align*}
			\lcm(z_1,w)&=\lcm(x_{h}(z/x_{i}),x_{h}(u/x_{i}))\\
			&=\lcm(x_{h}(x_{k}(u/x_{p})/x_{i}),x_{h}(u/x_{i}))\\
			&=x_{h}(x_{k}u)/x_{i}\in G(\HS_1(I)),
		\end{align*}
		and the property $(*)$ is satisfied as $h\in\Omega$, that is $\deg_{x_h}(u)<\deg_{x_h}(v)$ and as $h\ne k$, we have $\deg_{x_h}(u_1)=\deg_{x_h}(x_{k}u)=\deg_{x_h}(u)<\deg_{x_h}(v)\le\deg_{x_h}(x_{\ell}v)=\deg_{x_h}(v_1)$.
		
		\medskip\noindent
		\textsc{Subcase 2.2.2.} Suppose $z_2=x_{p}(z/x_{i})\in G(I)$. Note that 
		\begin{align*}
			z_2&=x_{p}(z/x_{i})=x_{p}(x_{k}(u/x_{p})/x_{i})=x_{k}(u/x_{i})
		\end{align*}
		and $d(z_2,w)=1$. Thus, Proposition \ref{prop:HS1(I)linquot} implies that
		\begin{align*}
			\lcm(z_2,w)&=\lcm(x_{k}(u/x_{i}),x_{h}(u/x_{i}))=x_{h}(x_{k}u)/x_{i}\in G(\HS_1(I)),
		\end{align*}
		and as before $\deg_{x_h}(u_1)<\deg_{x_h}(v_1)$. The proof is complete.
	\end{proof}
	
	By \cite[Theorem 2.4]{BanRam019}, a polymatroidal ideal $I\subset S$ has linear quotients with respect to the lex order $>_{\lex}$ (induced by any ordering of the variables $x_1,\dots,x_n$). Thus Corollary \ref{cor:HS1HSjLexEquig} implies that for all $j\ge0$,
	$$
	\HS_{j+1}(I)\ \subseteq\ \big(\HS_1(\HS_j(I))\big)_{>j+1}.
	$$
	
	Next we study when equality holds. This is the case when $I$ is actually matroidal, that is, a squarefree polymatroidal ideal \cite[Corollary 2.3]{Bay019}. The proof in \cite{Bay019} uses matroid and graph theory. We provide a purely algebraic proof. Firstly, we note the following general fact.
	\begin{Lemma}\label{Lem:HS>jSquarefree}
		Let $I\subset S$ be a squarefree ideal. Then, for all $j\ge0$,
		$$
		\big(\HS_1(\HS_j(I))\big)_{>j+1}\ =\ \HS_1(\HS_j(I)).
		$$
	\end{Lemma}
	\begin{proof}
		Since $I$ is squarefree, all $\HS_j(I)$ are squarefree. It follows from \cite[Lemma 4.4]{HHZ} that all monomials $w\in G(\HS_j(I))$ have $|\supp(w)|>j$.
		
		It is clear that $(\HS_1(\HS_j(I)))_{>j+1}\subseteq\HS_1(\HS_j(I))$. We show the opposite inclusion. Let $y\in G(\HS_1(\HS_j(I)))$. We claim that $|\supp(y)|>j+1$. By Proposition \ref{prop:HS1(I)linquot}, $y=\lcm(w_1,w_2)$ with $w_1,w_2\in G(\HS_j(I))$ such that $d(w_1,w_2)=1$. By Lemma \ref{Lem:d(u,v)=1}, $w_1=x_k(w_2/x_{\ell})$ for some $k\ne\ell$. Thus $y=\lcm(w_1,w_2)=x_\ell w_1$. We have shown that $w_1\in G(\HS_j(I))$ has $|\supp(w_1)|\ge j+1$. Since $\ell\notin\supp(w_1)$, $|\supp(y)|=1+|\supp(w_1)|\ge j+2>j+1$, as desired.
	\end{proof}
	
	\begin{Proposition}\label{prop:HSj+1=HS1HSjmatroidal}
		Let $I\subset S$ be a matroidal ideal. Then $\HS_{j+1}(I)=\HS_1(\HS_j(I))$ for all $j<\pd(I)$.
	\end{Proposition}
	\begin{proof}
		Since $I$ is squarefree, all homological shift ideals involved in the proof are squarefree. Fix the lex order $>_{\lex}$ induced by $x_1>\dots>x_n$. Then $I$ has linear quotients with respect to $>_{\lex}$, \cite[Theorem 2.4]{BanRam019}. For $u\in G(I)$, we denote by $\set(u)$ the following set $\{i:x_i\in (v\in G(I):v>_{\lex}u):u\}$. By equation (\ref{eq:HSi(I)linquot}),
		\begin{align}\label{eq:HSj(I)polymatroidal}
			\HS_j(I)\ =\ ({\bf x}_Fu\ :\ u\in G(I),\ F\subseteq\set(u),\ |F|=j).
		\end{align}
		For $j=0$, there is nothing to prove. Let $1\le j<\pd(I)$. By Proposition \ref{Prop:InclusionHSLinQuot} we have $\HS_{j+1}(I)\subseteq\HS_1(\HS_j(I))$. So we only need to prove that $\HS_1(\HS_j(I))\subseteq\HS_{j+1}(I)$. By Proposition \ref{prop:HS1(I)linquot}, we have
		$$
		\HS_1(\HS_j(I))=(\lcm(w_1,w_2)\ :\ w_1,w_2\in G(\HS_j(I)),\ d(w_1,w_2)=1).
		$$
		Thus, we must show that for all $w_1,w_2\in G(\HS_j(I))$ with $d(w_1,w_2)=1$ we have $\lcm(w_1,w_2)\in G(\HS_{j+1}(I))$. By equation (\ref{eq:HSj(I)polymatroidal}), $w_1={\bf x}_Fu\ne w_2={\bf x}_Gv$ with $u,v\in G(I)$, $F\subseteq\set(u)$, $G\subseteq\set(v)$ and $|F|=|G|=j$. Since $d(w_1,w_2)=1$, Lemma \ref{Lem:d(u,v)=1} gives $w_1=x_k(w_2/x_{\ell})$ for some $k\ne\ell$. As observed before, $w_1,w_2$ are squarefree. Hence, $\supp(w_1)=\{k\}\cup(\supp(w_2)\setminus\{\ell\})$. Note that as $\ell\in\supp(w_2)$, we can find $z\in G(I)$ such that $\ell\in\supp(z)$. Indeed if $\ell\in\supp(v)$ then we can choose $z=v$. Otherwise, $\ell\in G\subseteq\set(v)$, and $x_\ell v=x_sz$ with $z\in G(I)\setminus\{v\}$, $s\ne\ell$ and so $\ell\in\supp(z)$.
		
		Since $\ell\in\supp(z)\setminus\supp(u)$ it is $\deg_{x_\ell}(z)>\deg_{x_\ell}(u)$. By the dual exchange property $(**)$, the set of integers $h$ such that $\deg_{x_h}(u)>\deg_{x_h}(z)$ and $x_\ell(u/x_h)\in G(I)$ is non empty. Hence, the following set is non empty too
		$$
		\Omega\ =\ \{h\in[n]\setminus\{\ell\}\ :\ x_\ell(u/x_h)\in G(I)\}.
		$$
		
		\medskip\noindent
		\textsc{Case 1.} Assume there exists $h\in\Omega$ with $h>\ell$. Then $x_\ell(u/x_h)>_{\lex}u$ and $\ell\in\set(u)$. Now $\lcm(w_1,w_2)=\lcm(w_1,x_\ell(w_1/x_k))=x_{\ell}w_1=x_{\ell}{\bf x}_Fu$. Since $\ell\notin\supp(w_1)$, we also have $\ell\notin F$. Since $\ell\in\set(u)$, we have that $F\cup\{\ell\}$ is a subset of $\set(u)$ having cardinality $j+1$ and $\lcm(w_1,w_2)={\bf x}_{F\cup\{\ell\}}u\in G(\HS_{j+1}(I))$ by equation (\ref{eq:HSi(I)linquot}), as desired.
		
		\medskip\noindent
		\textsc{Case 2.} Assume that for all $h\in\Omega$ we have $h<\ell$. Choose some $h\in\Omega$. Then, $u>_{\lex}x_\ell(u/x_h)\in G(I)$. Set $w=x_\ell(u/x_h)$. Hence, this time $h\in\set(w)$. Note that $h\in\supp(u)$ and since $w_1={\bf x}_Fu$ is squarefree, $h\notin F$. We are going to show that $F\subseteq\set(w)$. Hence, we will have $F\cup\{h\}\subseteq\set(w)$ and then the desired conclusion:
		$$
		\lcm(w_1,w_2)=x_{\ell}{\bf x}_Fu=x_h{\bf x}_F(x_\ell(u/x_h))={\bf x}_{F\cup\{h\}}w\in G(\HS_{j+1}(I)).
		$$
		Let $m\in F$, then for some $p\ne m$, $x_m(u/x_p)\in G(I)$ and $x_m(u/x_p)>_{\lex}u$. So $m<p$.
		
		\medskip\noindent
		\textsc{Subcase 2.1.} Let $p=h$. Then $\ell>h=p>m$. Hence $x_m(u/x_p)=x_m(u/x_h)>_{\lex}x_\ell(u/x_h)=w$. Whence, $m\in\set(w)$ in this case, as desired.
		
		\medskip\noindent
		\textsc{Subcase 2.2.} Let $p\ne h$. Then $d(x_m(u/x_p),w)=d(x_m(u/x_p),x_\ell(u/x_h))=2$, $h\ne m$, $h\ne\ell$, $\ell\ne m$, $p\ne h$, $p\ne m$,  and
		\begin{align*}
			\deg_{x_h}(w)&<\deg_{x_h}(x_m(u/x_p)),& \deg_{x_\ell}(w)&>\deg_{x_\ell}(x_m(u/x_p)),\\
			\deg_{x_m}(w)&<\deg_{x_m}(x_m(u/x_p)),& \deg_{x_p}(w)&>\deg_{x_p}(x_m(u/x_p)).
		\end{align*}
		Whereas, for all $q\ne h,\ell,m,p$ we have $\deg_{x_q}(w)=\deg_{x_q}(x_m(u/x_p))$. Since $x_m(u/x_p)$, $w\in G(I)$ and $\deg_{x_m}(x_m(u/x_p))>\deg_{x_m}(w)$, by the dual exchange property $(**)$ we have either $x_m(w/x_\ell)\in G(I)$ or $x_m(w/x_p)\in G(I)$.
		
		\medskip\noindent
		\textsc{Subcase 2.2.1.} Assume $x_m(w/x_\ell)\in G(I)$. Note that
		$$
		x_m(w/x_\ell)=x_m(x_\ell(u/x_h))/x_\ell=x_m(u/x_h)\in G(I).
		$$
		Thus $m\in\Omega$ and by assumption $m<\ell$. Hence $x_m(w/x_\ell)>_{\lex}w$ and so $m\in\set(w)$.
		
		\medskip\noindent
		\textsc{Subcase 2.2.2.} Assume $x_m(w/x_p)\in G(I)$. In this case, since $m<p$ we have $x_m(w/x_p)>_{\lex}w$, and again $m\in\set(w)$. The proof is complete
	\end{proof}
	
	Proposition \ref{prop:HSj+1=HS1HSjmatroidal} and Theorem \ref{Thm:HS1(I)PolyMatr} yield another proof of Conjecture \ref{Conj:BBH} for matroidal ideals, one was already obtained in \cite[Theorem 2.2]{Bay019}.
	\begin{Corollary}\label{cor:HSj(I)matroidal}
		\textup{\cite[Theorem 2.2]{Bay019}} Let $I\subset S$ be a matroidal ideal. Then $\HS_j(I)$ is again a matroidal ideal for all $j$.
	\end{Corollary}
	
	Unfortunately, in general we could have $\HS_{j+1}(I)\ne(\HS_1(\HS_j(I)))_{>j+1}$.
	\begin{Example}\label{Ex:HS1(HSj)notTransversal}
		\rm Let $I=(x_1,x_2,x_3,x_4)(x_3,x_4,x_5)\subset S=K[x_1,\dots,x_5]$. Using the \textit{Macaulay2} \cite{GDS} package $\mathtt{HomologicalShiftIdeals}$ \cite{F2Pack} we verified that $\HS_j(I)$ is polymatroidal for all $j$. We have $\HS_{j+1}(J)\ne(\HS_1(\HS_j(J)))_{>j+1}$ for $j=1,2,3,4$. For instance $$(x_1x_2x_3x_4x_5)x_3^2x_4^2\in G(\HS_1(\HS_1(J)))\setminus G(\HS_2(J)).$$
	\end{Example}
	
	\textit{Acknowledgment.} I thank the referee whose suggestions greatly improved the readability and the quality of the paper.

\end{document}